\documentclass[final]{amsart}
\usepackage{amsmath,amsthm,amsfonts,amssymb,indentfirst,xspace}
\usepackage{xcolor}
\usepackage{graphicx,psfrag,subfigure}
\usepackage[notref,notcite]{showkeys}
\usepackage[ps2pdf]{hyperref}
\usepackage[initials,nobysame]{amsrefs}

\newenvironment{highlight}[1][yellow!70!black]{\color{#1}}{}

\newcommand{\slns}{\eqref{eNSXdef}--\eqref{eNSWebber}\xspace}
\newcommand{\nsequations}{\eqref{eNS}--\eqref{eDivFree}\xspace}
\newcommand{\mceuler}{\eqref{eFlowN}--\eqref{euN}\xspace}
\newcommand{\mceulert}{\eqref{eFlowit0}--\eqref{euNt0}\xspace}

\newcommand{\del}{\partial}
\newcommand{\lap}{\triangle}

\newcommand{\inv}{^{-1}}
\newcommand{\transpose}{^*}

\newcommand{\grad}{\nabla}
\newcommand{\gradt}{\grad\transpose}
\newcommand{\divergence}{\grad \cdot}
\newcommand{\curl}{\grad \times}

\renewcommand{\epsilon}{\varepsilon}
\renewcommand{\leq}{\leqslant}
\renewcommand{\geq}{\geqslant}

\newcommand{\E}{\mathbb{E}}
\newcommand{\F}{\mathcal{F}}
\DeclareMathOperator{\var}{Var}
\newcommand{\lhp}{\boldsymbol{\mathrm{P}}}
\newcommand{\as}{\text{a.s.}}

\newcommand{\W}{\mathcal{W}}

\newcommand{\R}{\mathbb{R}}

\newcommand{\N}{\mathbb{N}}
\newcommand{\holderspace}[2]{\ensuremath{C^{\ifx0#1{#2}\else{#1,#2}\fi}}}

\newcommand{\U}{\mathcal{U}}
\newcommand{\M}{\mathcal{M}}

\newif\iftextstyle
\textstyletrue
\everydisplay\expandafter{\the\everydisplay\textstylefalse}

\newcommand{\abs}[1]{\iftextstyle\lvert#1\rvert\else\left\lvert#1\right\rvert\fi}
\newcommand{\norm}[1]{\iftextstyle\lVert#1\rVert\else\left\lVert#1\right\rVert\fi}
\newcommand{\holdernorm}[3]{\norm{#1}_{\ifx0#2{#3}\else{#2,#3}\fi}}
\newcommand{\holdersnorm}[2]{\abs{#1}_{#2}}

\newcommand{\lpnorm}[2]{\norm{#1}_{\smash{L^{\!#2}_{\vphantom{h}}}\vphantom{L^{\!#2}}}}
\newcommand{\linfnorm}[1]{\lpnorm{#1}{\infty}}

\newcommand{\ip}[2]{\iftextstyle(\else\left(\fi #1, #2 \iftextstyle)\else\right)\fi}
\newcommand{\qv}[2]{\iftextstyle\langle\else\left\langle\fi #1, #2 \iftextstyle\rangle\else\right\rangle\fi}

\numberwithin{equation}{section}
\allowdisplaybreaks

\newtheorem{theorem}{Theorem}[section]
\newtheorem{lemma}[theorem]{Lemma}
\newtheorem{proposition}[theorem]{Proposition}
\newtheorem{corollary}[theorem]{Corollary}

\newtheorem*{theorem*}{Theorem}
\newtheorem*{lemma*}{Lemma}
\newtheorem*{proposition*}{Proposition}
\newtheorem*{corollary*}{Corollary}

\theoremstyle{definition}
\newtheorem{definition}[theorem]{Definition}

\theoremstyle{remark}

\newtheorem*{remark*}{Remark}

\begin{document}
\title[A stochastic-Lagrangian particle system for Navier-Stokes]{A
  stochastic-Lagrangian particle system for the Navier-Stokes
  equations.}  \author{Gautam Iyer} \address{%
  Department of Mathematics\\
  Stanford University} \email{gi1242@stanford.edu}  \author{Jonathan Mattingly} \address{%
  Department of Mathematics\\
  Duke University} \email{jonm@math.duke.edu}\thanks{GI was
  partly supported by NSF grant (DMS-0707920), and thanks the
  mathematics department at Duke for its hospitality. JCM was partial
  supported by an NSF PECASE awawrd (DMS-0449910) and a Sloan
  foundation fellowship.} 
\begin{abstract}
  This paper is based on a formulation of the Navier-Stokes equations
  developed by P. Constantin and the first author
  (\texttt{arxiv:math.PR/0511067}, to appear), where the velocity
  field of a viscous incompressible fluid is written as the expected
  value of a stochastic process. In this paper, we take $N$ copies of
  the above process (each based on independent Wiener processes), and
  replace the expected value with $\frac{1}{N}$ times the sum over
  these $N$ copies. (We remark that our formulation requires one to
  keep track of $N$ stochastic flows of diffeomorphisms, and not just
  the motion of $N$ particles.)

  We prove that in two dimensions, this system of interacting
  diffeomorphisms has (time) global solutions with initial data in the
  space $\holderspace{1}{\alpha}$ which consists of differentiable
  functions whose first derivative is $\alpha$ H\"older continuous
  (see Section \ref{sGexist} for the precise definition). Further, we
  show that as $N \to \infty$ the system converges to the solution of
  Navier-Stokes equations on any finite interval $[0,T]$. However for
  fixed $N$, we prove that this system retains roughly
  $O(\frac{1}{N})$ times its original energy as $t \to \infty$. Hence
  the limit $N \rightarrow \infty$ and $T\rightarrow \infty$ do not
  commute. For general flows, we only provide a lower bound to this
  effect. In the special case of shear flows, we compute the behaviour
  as $t \to \infty$ explicitly.
\end{abstract}
\maketitle
\section{Introduction}
The Navier-Stokes equations
\begin{gather}
\label{eNS} \del_t u + (u \cdot \grad) u - \nu \lap u + \grad p = 0\\
\label{eDivFree} \divergence u = 0
\end{gather}
describe the evolution of a velocity field of an incompressible fluid
with kinematic viscosity $\nu > 0$. These equations have been used to
model numerous physical problems, for example air flow around an
airplane wing, ocean currents and meteorological phenomena to name a
few \cites{bPedlosky,bChorinMarsden,bLamb}. The mathematical theory
(existence and regularity \cites{bibConstBook,bibLadyzhenskaya}) of
these equations have been extensively studied and is still one of the
outstanding open problems in modern PDE's \cites{constOP,clay}.

The questions addressed in this paper are motivated by a formalism of
\nsequations developed in \cite{detsns} (equations \slns below). This
formalism essentially superimposes Brownian motion onto particle
trajectories, and then averages with respect to the Wiener measure. In
this paper, we take $N$ independent copies of the Wiener process and
replace the expected value in the above formalism with $\frac{1}{N}$
times the sum over these $N$ independent copies (see equations
\mceuler for the exact details). In the original formulation the
random trajectory of a particle induced by a single Brownian motion interacts
with its own law. This is essentially a self-consistent, mean-field
interaction. In this paper, we replace this with $N$ copies or replica
whose average is used to approximate the interaction with the
processes own law.
This technique has been extensively used in numerical computation
(e.g. \cites{montecarlo,bRobertCasella,bibBremaud}). We remark
that in our formulation, we are required to keep track of
$N$ stochastic flows of diffeomorphisms, and not just the motion of $N$
different particles, as is the conventional approach.

We study both the behaviour as $N \to \infty$ and $t \to \infty$ of
the system obtained. The behaviour as $N \to \infty$ is as expected:
In two dimensions on any finite time interval $[0,T]$, the system
converges as $N \rightarrow \infty$ to the solution of the true
Navier-Stokes equations at rate roughly $O(\frac{1}{\sqrt{N}})$. In
three dimensions, we can only guarantee this if we have certain
apriori bounds on the solution (Theorem \ref{tConvergence}). These
apriori bounds are of course guaranteed for short time, but are
unknown (in the $3$-dimensional setting) for long time
\cites{clay,constOP}.

\begin{figure}[thb]
\subfigure[$N=2$]{
    \includegraphics[width=5.9cm]{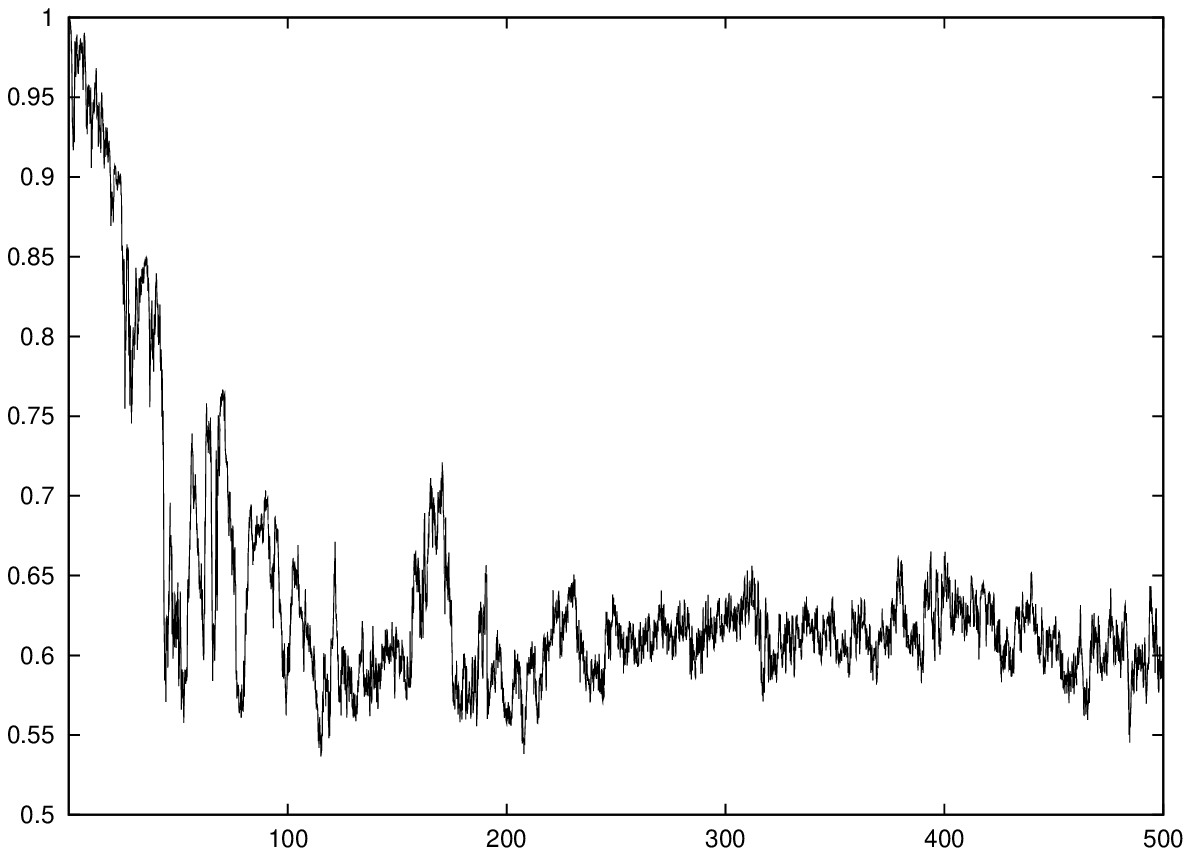}\label{fn2}
  }\quad \subfigure[$N=8$]{
    \includegraphics[width=5.9cm]{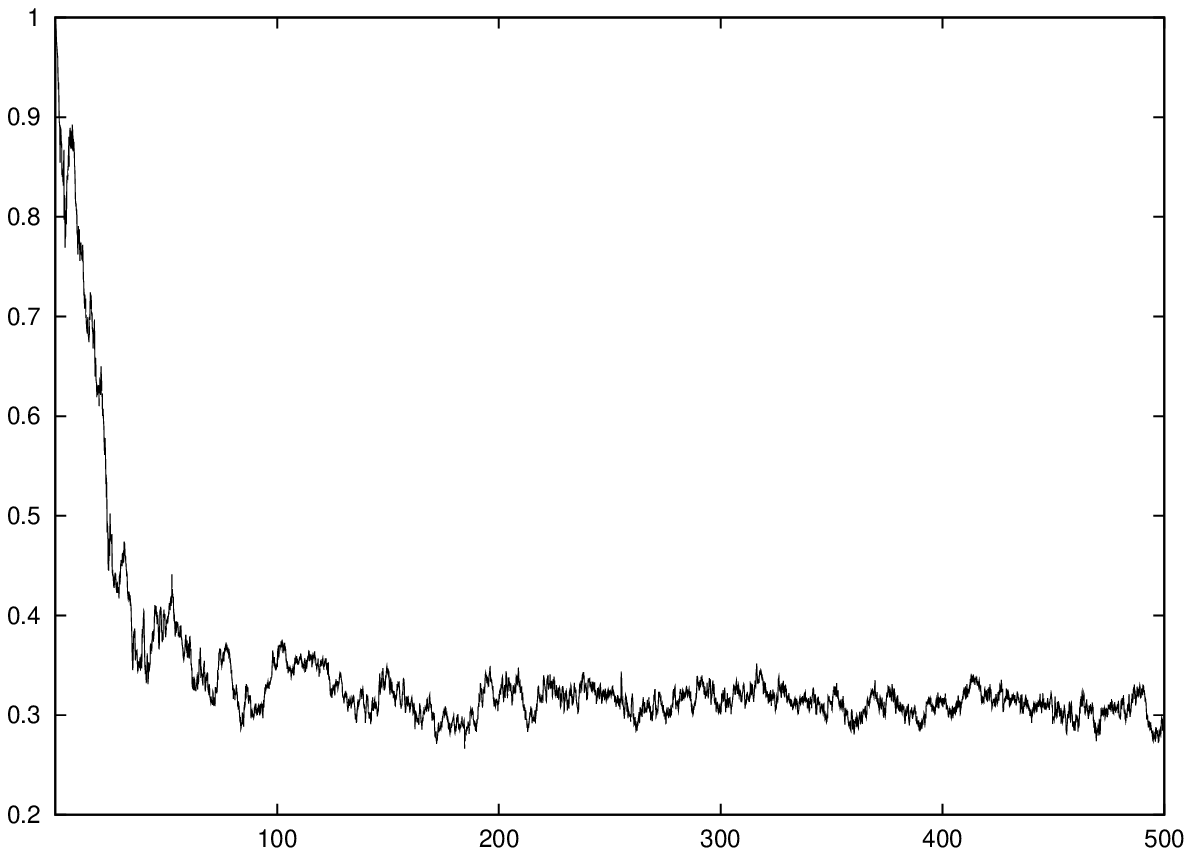}\label{fn8}
}
\caption{Graph of $\norm{\omega^N_t}_{L^2}^2$ vs time}
\end{figure}

At first glance, the behaviour as $t \to \infty$ for fixed $N$ is less
intuitive. For the $2$-dimensional problem, Figures \ref{fn2} and
\ref{fn8} show a graph of $\lpnorm{\omega^N_t}{2}^2$ vs time, with $N
= 2$ and $N = 8$ respectively\footnote{These computations were done using a $24 \times 24$ mesh on the periodic box with side length $2\pi$. The initial vorticity was randomly chosen, and normalized with $\lpnorm{\omega_0}{2} = 1$. The behaviour depicted in these two figures is however characteristic, and insensitive to changes of the mesh size, length, or diffusion coefficient}. A little reflection shows that this
behavior is not completely surprising. The dissipation occurs through
the averaging of different copies of the flow. With only $N$ copies,
one can only produce dissipation of order $1/N$ of the original
energy. It is tempting, to speak of the inability to represent the
correct interaction of small scale structures with such small number
of data. However, we can not make this precise and since each of the
objects being averaged is an entire diffeomorphism with an infinite
amount of information it is unclear what this means.

In Section \ref{sTtoInfty} we obtain a sharp lower bound to this effect. We show (Theorem
\ref{tLimSupGradULowerBound}) that
$$
\limsup_{t \to \infty} \E \lpnorm{\grad u_t}{2}^2 \geq \frac{1}{N L^2} \lpnorm{u_0}{2}^2,
$$
where $L$ is a length scale.
Further, we explicitly compute the $t \to \infty$ behaviour in the special case of shear flows and verify that our lower bound is sharp.\medskip

We remark that we considered the analogue of the system above
for the one dimensional Burgers equations. As is well known the
viscous Burgers equations have global strong solutions. However
preliminary numerical simulations show that the system above forms
shocks almost surely, even for very large $N$. We are currently
working on understanding how to continue this system past these
shocks, in a manner analogous to the entropy solutions for the
inviscid Burgers equations, and studying its behaviour as $t \to
\infty$ and $N \to \infty$.

We do not propose this particle system as an efficient particle method
for numerical computation. Though there may be special cases were it
may be useful, in general the computational cost of representing $N$
entire diffeomorphisms is large. Rather we see it as an interesting
and novel regularization which might give useful insight into the
structure and role of dissipation in the system. 
\section{The particle system}\label{sParticleSystem}

In this section we construct a particle system for the Navier-Stokes
equations based on stochastic Lagrangian trajectories. We begin by
describing a stochastic Lagrangian formulation of the Navier-Stokes
equations developed in \cites{thesis,detsns}.

Let $W$ be a standard $2$ or $3$-dimensional Brownian motion, and
$u_0$ some given divergence free $\holderspace{2}{\alpha}$ initial
data. Let $\E$ denote the expected value with respect to the Wiener
measure and $\lhp$ be the Leray-Hodge projection onto divergence free
vector fields. Consider the system of equations
\begin{align}
  \label{eNSXdef} dX_t(x) &= u_t ( X_t(x) ) \, dt + \sqrt{2\nu}\,dW_t, \qquad X_0(x) = x,\\
  \label{eNSWebber} u_t &= \E\, \lhp\left[ (\gradt Y_t) (u_0 \circ Y_t)
  \right],\qquad Y_t = X_t\inv.
\end{align}
With a slight abuse of notation, we denote by $X_t$ the map from
initial conditions to the value at time $t$. Hence $X_t$ is a
stochastic flow of diffeomorphisms with $X_0$ equal to the identity
and $Y_t$ the ``spatial'' inverse. In other words, $Y_t\colon
X_t(x)\mapsto x$. Also by $\gradt Y_t$ we mean the transpose of the
Jacobian of map $Y_t$.  Observe that $(\gradt Y_t) (u_0 \circ Y_t)$
can be viewed as a function of $x$ where both the Jacobian and the
vector field $u_0 \circ Y_t$ to which it is applied are both evaluated
at $x$.  We impose periodic boundary conditions on the displacement
$\lambda_t(y) = X_t(y) - y$, and on the Leray-Hodge projection $\lhp$.

In \cites{thesis,detsns} it was shown that the system \slns is
equivalent to the Navier-Stokes equations in the following sense: If
the initial data is regular ($\holderspace{2}{\alpha}$), then the pair
$X,u$ is a solution to the system \slns if and only if $u$ is a
(classical) solution to the incompressible Navier-Stokes equations
with periodic boundary conditions and initial data $u_0$.

We digress briefly and comment on the physical significance of
\slns. Note first that equation \eqref{eNSWebber} is
\textit{algebraically} equivalent to the equations
\begin{gather}
\label{eBiotSavart} u_t = - \lap\inv \curl \omega_t\\
\label{eVorticityTransport} \omega_t = \E \left[ (\grad X_t) \omega_0 \right] \circ Y_t
\end{gather}
This follows by direct computation, and was shown \cite{ele} and
\cite{detsns} for instance. We recall that \eqref{eVorticityTransport}
is the usual vorticity transport equation for the Euler equations, and
\eqref{eBiotSavart} is just the Biot-Savart formula.

Thus in particular when $\nu = 0$, the system \slns is exactly the
incompressible Euler equations. Hence the system \slns essentially
does the following: We add Brownian motion to Lagrangian
trajectories. Then recover the velocity $u$ in the same manner as for
the Euler equations, but additionally average out the noise.

We remark that the system \slns is non-linear in the sense of McKean
\cite{bibSznitman}. The drift of the flow $X$ depends on its
distribution. However in this case, the law of $X$ alone is not enough
to compute the drift $u$. This is because of the presence of the
$\gradt Y$ term in \eqref{eNSWebber}, which requires knowledge of
spatial covariances, in addition to the law of $X$. In other words,
one needs that law of the entire flow of diffeomorphism and not just
the law of the one-point motions.  \medskip

We now motivate our particle system. For the formulation \slns above,
the natural numerical scheme would be to use the law of large numbers
to compute the expected value. Let $(W^i)$ be a sequence of
independent Wiener processes, and consider the system
\begin{align}
  \label{eFlowN} dX^{i,N}_t &= u^N_t \left( X^{i,N}_t \right) \, dt + \sqrt{2\nu}\,dW^i_t,\quad  Y^{i,N}_t = \left( X^{i,N}_t \right)\inv\\
\label{euN} u^N_t &= \frac{1}{N} \sum_{i=1}^N \lhp\left[ (\gradt Y^{i,N}_t) (u_0 \circ Y^{i,N}_t) \right]
\end{align}
with initial data $X_0(x) = x$. We impose again periodic boundary
conditions on the initial data $u_0$, the displacement $\lambda_t(x) =
X_t(x) - x$, and the Leray-Hodge projection $\lhp$.

We remark that the algebraic equivalence of \eqref{eNSWebber} and
\eqref{eBiotSavart}--\eqref{eVorticityTransport} is still valid in
this setting. Thus the system \mceuler could equivalently be
formulated by replacing equation \eqref{euN} with the more familiar
equations
\begin{equation}
\label{eOmegaN} \omega_t^N = \frac{1}{N} \sum_{i = 1}^N \left[ (\grad X^{i,N}_t) \omega_0\right] \circ Y^{i,N}_t
\end{equation}
and
\begin{equation}
  \label{eBiotSavartN} u^N_t = - \lap\inv \curl \omega^N_t.
\end{equation}

Finally we clarify our previous remark, stating that the above formulation requires us to keep track of $N$ stochastic flows, and the knowledge of the one point motions of $X^{i,N}_t$ alone is not sufficient. The standard method of obtaining a solution to the heat equation (assuming the drift $u$ is time independent) would be to consider the process \eqref{eNSXdef}, and read off the solution $\theta$ by
\begin{equation*}
\theta_t(a) = \E \theta_0 \circ X_t(a)
\end{equation*}
where $\theta_0$ is the given initial temperature distribution.
Thus knowing the trajectories (and distribution) of the process $X$ starting at one particular point $a$ will be sufficient to determine the solution $\theta_t$ at that point.

This however is not the case for our representation. The reason is twofold: First, the representation \slns involves a non-local singular integral operator. Second our representation involves composing with the \textit{spatial inverse} of the flow $X_t$, and then averaging. If we for a moment ignore the non-locality of the Leray-Hodge projection, determining $u_t$ at one fixed point $a$ one would need the law of $Y_t(a)$, for which the knowledge of $X_t(a)$ alone is not enough. One needs the entire (spatial) map $X_t$ to compute the spatial inverse $Y_t(a)$.

The above is not a serious impediment to a numerical implementation. Given an initial mesh $\Delta$, we first compute $X^{i,N}_t$ on this mesh. By definition of the inverse, one knows $Y^{i,N}_t$ on the (non-uniform) mesh $X^{i,N}_t(\Delta)$, after which one can interpolate and find $Y^{i,N}_t$ on the mesh $\Delta$. In two spatial dimensions, global existence and regularity (Theorem \ref{t2DGlobalExistence}) together with incompressibility will show that this mesh does not degenerate in finite time.

This surprisingly is \textit{not} the case for the (one dimensional)
Burgers equations. Numerical computations indicate the mesh $\Delta$
almost surely degenerates in finite time for non-monotone initial
data, and the solution `shocks' almost surely. Thus, while the system
\mceuler appears natural, and convergence as $N \to \infty$ is to be
expected, caution is to be exercised. We suspect that the results
(existence, convergence, etc.) proven in this paper for the system
\mceuler are in fact \textit{false} for the Burgers equations. This is
indeed puzzling as global existence, and regularity for the viscous
Burgers equations is well known. It further underlines the fact that
the finite $N$ approximation modifies the dissipation in a different
way then other approximation such as a spectral
approximation. \medskip

In the next section, we show the existence of global solutions to
\mceuler in two dimensions. In section \ref{sNtoInfty}, we show that
the solution to \mceuler converges to the solution of the
Navier-Stokes equations as $N \to \infty$. Finally in Section
\ref{sTtoInfty} we study the behaviour of the system \mceuler as $t
\to \infty$ (for fixed $N$), and partially explain the behaviour shown
in Figures \ref{fn2} and \ref{fn8}.

\section{Global existence of the particle system in two dimensions.}\label{sGexist}
In this section we prove that the particle system \mceuler has global
solutions in two dimensions. Once we are guaranteed global in time
solutions, we are able to study the behaviour as $t \to \infty$, which
we do in Section \ref{sTtoInfty}. We remark also that as a consequence
of Theorem \ref{t2DGlobalExistence} (proved here), our convergence
result as $N \to \infty$ (Theorem \ref{tConvergence}) applies on any
finite time interval $[0,T]$ in the two dimensional situation.\medskip

We first establish some notational convention: We let $L > 0$ be a
length scale, and assume work with the spatial domain is $[0, L]^d$,
where $d\geq 2$ is the spatial dimension. We define the
non-dimensional $L^p$ and H\"older norms by
\begin{gather*}
  \lpnorm{u}{p} = \left( \frac{1}{L^d} \int_{[0,L]^d} \abs{u}^p \right)^{\frac{1}{p}}\\
  \lpnorm{u}{\infty} = \sup_{x \in [0,L]^d} \abs{u (x) }\\
  \holdersnorm{u}{\alpha} = L^\alpha \sup_{x\neq y \in [0,L]^d} \frac{ \abs{u(x) - u(y)} }{ \abs{x - y}^\alpha }\\
  \holdernorm{u}{k}{\alpha} = \sum_{\abs{m} = k} L^k \holdersnorm{D^m
    u}{\alpha} + \sum_{\abs{m} < k} L^{\abs{m}} \lpnorm{D^m u}{\infty}
\end{gather*}
We denote the H\"older space $\holderspace{k}{\alpha}$ and Lebesgue
space $L^p$ to be the space of functions $u$ which are periodic on
$[0,L]^d$ with $\holdernorm{u}{k}{\alpha} < \infty$ and $\lpnorm{u}{p}
< \infty$ respectively.

Let $W^1, \dots, W^N$ be $N$ independent ($2$ dimensional) Wiener
processes, with filtration $\F_t$.
\begin{proposition}\label{pLocalExistence}
  Let $t_0 \geq 0$, and $u_{t_0}^1, \dots, u_{t_0}^N$ be
  $\F_{t_0}$-measurable, periodic mean $0$ functions such that the
  norms $\holdernorm{u^i_{t_0}}{1}{\alpha}$ are almost surely
  bounded. Then there exist $T = T( \alpha,
  \holdernorm{u_{t_0}^1}{1}{\alpha}, \dots,
  \holdernorm{u_{t_0}^N}{1}{\alpha} )$ such that the system
\begin{align}
  \label{eFlowit0} X^i_{t_0, t}(x) &= x + \int_{t_0}^t u_s( X^i_{t_0,s}(x) )
  \, ds + \sqrt{2\nu} (W^i_t - W^i_{t_0}), \quad   Y^i_{t_0, t} = (X^i_{t_0,t})\inv\\
\label{euNt0} u_t(x) &= \frac{1}{N} \sum_{i=1}^N \lhp\left[ (\gradt Y^i_{t_0, t})
  \, (u^i_{t_0} \circ Y^i_{t_0, t}) \right](x) 
\end{align}
has a solution on the interval $[t_0, T]$. Further there exists deterministic $U =
U(\alpha, \holdernorm{u^i_{t_0}}{1}{\alpha})$ such that
\begin{equation}\label{eUiBound}
  \sup_{t \in [t_0, T]} \holdernorm{\lhp\left[ (\gradt Y^i_{t_0, t})
      \, u^i_{t_0} \circ Y^i_{t_0, t} \right]}{1}{\alpha} \leq U 
\end{equation}
almost surely. Consequently, $u \in C( [t_0, T], \holderspace{1}{\alpha} )$, and $\holdernorm{u}{1}{\alpha} \leq U$.
\end{proposition}
Proposition \ref{pLocalExistence} is proved in Appendix
\ref{aLocalExistence}

\begin{definition}
  We call $X^i$, $u$ in Proposition \ref{pLocalExistence} the solution
  of the system $\mceulert$ with initial data $u^1_{t_0}, \dots,
  u^N_{t_0}$.
\end{definition}

The functions $X^i, Y^i$ are not periodic themselves, but have
periodic displacements: Namely, if we define
\begin{gather}
\label{eLambdai} \lambda^i_t(y) = X^i_t(y) - y\\
\label{eMui} \mu^i_t(x) = Y^i_t(x) - x\;,
\end{gather}
then $\mu^i, \lambda^i$ are periodic.

We remark that if $t_0 = 0$ and the $\omega^i_{t_0}$'s are all equal,
then the system \mceulert reduces to \mceuler. However when formulated
as above, solutions can be continued past time $T$ by restarting the
flows $X^i$, as in the Lemma below.

\begin{lemma}\label{lMarkov}
  Say $t_0 < t_1 < t_2$, and $X^i_{t_0,s}$, $u_s$ solve \mceulert on
  $[t_0, t_2]$, with initial data $u^i_{t_0}$. For any $t > t_0$
  define
$$
u^i_t = \lhp \left[ (\gradt Y^i_{t_0, t}) \, u^i_{t_0} \circ
  Y^i_{t_0,t} \right].
$$
Let $\tilde X^i_{t_1,s}$, $\tilde u_s$ solve \mceulert on $[t_1, t_2]$
with initial data $u^i_{t_1}$. Then for all $s \in [t_1, t_2]$ we have
$\tilde u_s^i = u^i_s$ and
$$
X^i_{t_0,s} = \tilde X^i_{t_1, s} \circ X^i_{t_0, t_1}
$$
almost surely.
\end{lemma}

\begin{proof}
  The proof is identical to the proof of Proposition 3.3.1 in
  \cite{thesis}, and we do not provide it here.
\end{proof}

For the remainder of this section, we assume without loss of
generality that $t_0 = 0$ (we allow of course $\F_0$ to be
non-trivial). For ease of notation, we use $X_s$ to denote
$X_{0,s}$. We now prove that the system \mceulert has global solutions
in two dimensions. This essentially follows from a Beale-Kato-Majda
type condition \cite{bkm}, and the two dimensional vorticity
transport.

\begin{lemma}\label{lSIOlogBound}
  If $u$ is divergence free and periodic in $\R^d$, then for any
  $\alpha \in (0,1)$, there exists a constant $c = c(\alpha, d)$ such
  that
$$
\linfnorm{\grad u} \leq c \linfnorm{\omega} \left( 1 + \ln^+ \left(
    \frac{ \holdersnorm{\omega}{\alpha} }{ \linfnorm{\omega} } \right)
\right)
$$
where $\omega = \curl u$.
\end{lemma}

The lemma is a standard fact about singular integral operators, and we
provide a proof in Appendix \ref{aSIOlogBound} for completeness.

\begin{theorem}\label{t2DGlobalExistence}
  In $2$ dimensions, the system \mceulert has time global solutions
  provided the initial data $u^1_0, \dots u^N_0$ is periodic,
  $\F_0$-measurable with $\holdernorm{u^i_0}{1}{\alpha}$ bounded
  almost surely. In particular we have global solutions to \mceulert
  in two dimensions, if $u^1_0 = \cdots = u^N_0 = u_0$ is
  deterministic, H\"older $1 + \alpha$ and periodic.
\end{theorem}
\begin{proof}
Taking the curl of \eqref{euNt0} gives the familiar Cauchy formula \cites{ele,detsns,sperturb}
\begin{equation}\label{eCauchy3D}
  \omega_t = \frac{1}{N} \sum_{i=1}^N \left[ (\grad X^i_t) \, \omega^i_0 \right] \circ Y^i_t,
\end{equation}
where $\omega_t = \curl u_t$. In two dimensions reduces to
\begin{equation}\label{eCauchy2D}
\omega_t = \frac{1}{N} \sum_{i=1}^N \omega^i_0 \circ Y^i_t.
\end{equation}
Taking H\"older norms gives
\begin{equation}\label{eOmegaHolderBound}
  \holdernorm{\omega_t}{0}{\alpha} \leq \frac{c}{N} \sum_{i=1}^N
  \holdernorm{\omega^i_0}{0}{\alpha} (1 + \linfnorm{\grad Y^i}^\alpha)
  \quad\as 
\end{equation}

Now differentiating \eqref{eFlowit0} gives
$$
\grad X^i_t = I + \int_0^t (\grad u_s \circ X^i_s) \, \grad X^i_s \, ds \quad\as
$$
Taking the $L^\infty$ norm, and applying Gronwall's Lemma shows
$$
\linfnorm{\grad X^i_t} \leq \exp\left( c \int_0^t \linfnorm{ \grad u_s
  } \, ds \right) \quad\as 
$$
Recall $\divergence u = 0$, and hence $\det(\grad X^i) = 1$ almost
surely. Thus the entries of $\grad Y$ are a polynomial (of degree $1$)
in the entries of $\grad X$. This immediately gives
\begin{equation}\label{eGradYiLinfBound}
\linfnorm{\grad Y^i_t} \leq \exp\left( c \int_0^t \linfnorm{ \grad u_s } \, ds \right)
\end{equation}
almost surely. Combining this with \eqref{eOmegaHolderBound} gives us the apriori bound
\begin{equation}\label{eOmegaGradUBound}
  \holdernorm{\omega_t}{0}{\alpha} \leq \frac{c}{N} \sum_{i=1}^N
  \holdernorm{\omega^i_0}{0}{\alpha} \exp\left( c \int_0^t \linfnorm{
      \grad u_s } \, ds \right) 
\end{equation}

Applying Lemma \ref{lSIOlogBound} gives us
\begin{align*}
  \linfnorm{ \grad u } &\leq c \linfnorm{\omega} \left( 1 +
    \ln^+\left( \frac{\holdersnorm{\omega}{\alpha}
      }{\linfnorm{\omega}}\right) \right)\\ 
  &\leq c \linfnorm{\omega} \left( 1 + \ln^+\left(
      \frac{\holdersnorm{\omega}{\alpha} }{\frac{1}{N} \sum_{i=1}^N
        \linfnorm{\omega^i_0}}\right) + \ln^+\left(
      \frac{\frac{1}{N}\sum_{i=1}^N \linfnorm{\omega^i_0}
      }{\linfnorm{\omega}}\right) \right)
\end{align*}
Note that the function $x \ln^+ \frac{1}{x}$ is bounded, so the last
term on the right can be bounded above by some constant $c_0$. For the
remainder of the proof, we let $c_0 = c_0(
\holdernorm{\omega^i_0}{0}{\alpha}, \alpha)$ denote a constant (with
dimensions that of $\omega$) which changes from line to line. Thus
\begin{align*}
  \linfnorm{\grad u_t} &\leq c_0 + c \linfnorm{\omega_t} \left( 1 +
    \ln^+\left( \frac{\holdersnorm{\omega_t}{\alpha} }{\frac{1}{N}
        \sum_{i=1}^N \linfnorm{\omega^i_0}}\right) \right)\\ 
  &\leq c_0 + c \linfnorm{\omega_t} \left( 1 + \int_0^t
    \linfnorm{\grad u_s} \, ds \right)
\end{align*}
and hence
$$
\linfnorm{\grad u_t} \left( 1 + \int_0^t \linfnorm{\grad u_s} \, ds
\right)\inv \leq c_0 + c \linfnorm{\omega_t}.
$$
Integrating gives us the apriori bound
\begin{align}
  \nonumber \int_0^t \linfnorm{\grad u_t} &\leq \exp\left(c_0 t + c
    \int_0^t \linfnorm{\omega_t} \right) - 1\\ 
\label{eAprioriIntGradULinf}    &\leq c_0 t e^{c_0 t}
\end{align}
since \eqref{eCauchy2D} implies $\linfnorm{\omega_t} \leq \frac{1}{N}
\sum \linfnorm{\omega^i_0}$. 

Now if we set $\omega^i_t = \omega_0^i \circ Y^i_t$, then
\eqref{eGradYiLinfBound} and the apriori bound
\eqref{eAprioriIntGradULinf} gives
\begin{equation}\label{eOmegaIAprioriHolder}
\holdernorm{\omega^i_t}{0}{\alpha} \leq \holdernorm{\omega^i_0}{0}{\alpha} \exp\left( c_0 t e^{c_0 t} \right)
\end{equation}
If $u^i_t$ is as in Lemma \ref{lMarkov}, then \eqref{eOmegaIAprioriHolder} shows
\begin{equation}\label{eGradUBound}
  \holdernorm{\grad u^i_t}{0}{\alpha} \leq c
  \holdernorm{\omega^i_t}{0}{\alpha} \leq
  c\holdernorm{\omega_0}{0}{\alpha} \exp\left( c_0 t e^{c_0 t}
  \right). 
\end{equation}
Since the mean velocity is a conserved quantity, a bound on
$\holdernorm{\grad u^i_t}{0}{\alpha}$ immediately gives a bound on
$\holdernorm{u^i_t}{1}{\alpha}$, which in conjunction with local
existence (Proposition \ref{pLocalExistence}), and Lemma \ref{lMarkov}
concludes the proof.
\end{proof}

\section{Convergence as \texorpdfstring{$N\to\infty$}{N to infinity}}\label{sNtoInfty}
In this section, we fix a time interval $[0,T]$, and show that the
particle system \mceuler converges to the solution to the
Navier-Stokes equations as $N \to \infty$. The rate of convergence is
$O(\frac{1}{\sqrt{N}})$, which is comparable to the convergence rate of the random vortex method \cites{bibLong,bibChorin}. As mentioned earlier, the system is intrinsically non-local, and propagation of chaos \cite{bibSznitman} type estimates are not easy to obtain. Consequently convergence results based on spatially averaged norms are easier to obtain, and we present one such result in this section, under assumptions which are immediately guaranteed by local existence.

\begin{theorem}\label{tConvergence}
  For each $i, N$, let $X^{i,N}$, $u^N$ be a solution to the particle
  system \mceuler with initial data $u_0$ on some time interval
  $[0,T]$. Let $u$ be a solution to the Navier-Stokes equations (with
  the same initial data) on the interval $[0,T]$. Suppose $U$ is such
  that
$$
\sup_{t \in [0,T]} \lpnorm{\grad u_t}{2} \leq \frac{U}{L}
\quad\text{and}\quad \sup_{t \in [0,T]} \lpnorm{ \grad u^{i,N}_t }{2}
\leq \frac{U}{L} \quad\as
$$
Then $(u^N) \to u$ in the following sense:
$$
\lim_{N \to \infty} \sup_{t \in [0,T]} \E \lpnorm{u^N_t - u_t}{2} = 0
$$
\end{theorem}

We remark that given $\holderspace{1}{\alpha}$ initial data, local
existence (Proposition \ref{pLocalExistence}) guarantees that the
conditions of this theorem are satisfied on some small interval
$[0,T]$. In two dimensions, Theorem \ref{t2DGlobalExistence} shows
that the conditions of this theorem are satisfied on any interval
finite $[0,T]$.

The proof will follow almost immediately from the following Lemma.
\begin{lemma}\label{lSPDEui}
Let $u^{i,N}_t = \lhp[ (\gradt Y^{i,N}_t) u_0\circ Y^{i,N}_t]$ be the $i^\text{th}$ summand in \eqref{euN}. Then $u^{i,N}$ satisfies the SPDE
\begin{multline}\label{eSPDEuiN}
  d u^{i,N}_t + \left[ (u^N_t \cdot \grad) u^{i,N}_t - \nu \lap
    u^{i,N}_t + (\gradt u^N_t) u^{i,N}_t +  \grad p^{i,N}_t\right] \,
  dt +\\ 
  + \sqrt{2\nu} \grad u^{i,N}_t dW^i_t = 0 
\end{multline}
and $u^N$ satisfies the SPDE
\begin{equation}\label{eSPDEuN}
  d u^N_t + \left[ (u^N_t \cdot \grad) u^N_t - \nu \lap u^N_t + \grad
    p^N_t \right] \, dt + \frac{\sqrt{2\nu}}{N} \sum_{i=1}^N \grad
  u^{i,N}_t dW^i_t = 0 
\end{equation}
\end{lemma}
\begin{remark*}
  We draw attention to the fact that the pressure term in
  \eqref{eSPDEuN} has bounded variation in time.
\end{remark*}
\begin{proof}
  We first recall a fact from \cites{detsns,thesis} (see also
  \cite{krylov-rosovski}). If $X$ is the stochastic flow
$$
dX_t = u_t\, dt + \sqrt{2\nu} \,dW_t
$$
and $Y_t = X_t \inv$ is the spatial inverse. Then the process $\theta_t
= f(Y_t)$ satisfies the SPDE
\begin{equation}\label{eSPDETheta}
  d\theta_t + (u_t \cdot \grad) \theta_t \, dt - \nu \lap \theta_t \,
  dt + \sqrt{2 \nu} \grad \theta_t \, dW_t = 0 
\end{equation}

This immediately shows that $Y^{i,N}_t$ and $v^{i,N}_t = u_0 \circ
Y^{i,N}_t$ both satisfy the SPDE \eqref{eSPDETheta}. For notational
convenience, we momentarily drop the $N$ as a superscript and use the
notation $v^{i,j}$ to denote the $j^\text{th}$ component of $v^i$.

Now we set $w^i_t = (\gradt Y^i_t) v^i_t$ and apply It\^o's formula:
\begin{align*}
  dw^{i,j}_t &= d(\del_j Y^i_t) \cdot v^i_t + (\del_j Y^i_t)\cdot dv^i_t +
  d\qv{\del_j Y^{i,k}_t}{v^{i,k}_t}\\ 
  &= (\del_j Y^i_t) \cdot \left[ -(u_t \cdot \grad) v^i_t + \nu \lap v^i_t
  \right] \, dt - \sqrt{2\nu} \del_j Y^i_t \cdot \left( \grad v^i_t \,
    dW^i_t\right) +\\ 
  &\qquad + v^i_t \cdot \left[ -((\del_j u_t) \cdot \grad) Y^i_t - (u_t \cdot
    \grad) \del_j Y^i_t + \nu \lap \del_j Y^i_t \right] \, dt -\\ 
  &\qquad - \sqrt{2\nu} v^i_t \cdot \left( \grad \del_j Y^i_t \, dW^i_t
  \right) + 2\nu \del^2_{jl} Y^{i,k}_t \del_l v^{i,k}_t \, dt\\ 
  &= \left[ -(u_t\cdot \grad)w^i_t + \nu \lap w^i_t - (\gradt u_t) \cdot w^i_t
  \right]\,dt - \sqrt{2\nu} \grad w^i_t\,dW^i_t
\end{align*}
Restoring the dependence on $N$ to our notation,
since $u^{i,N} = \lhp w^{i,N}$ we know that
$$
u^{i,N}_t = w^{i,N}_t + \grad q^{i,N}_t
$$
for some function $q^{i,N}_t$. Thus
\begin{align*}
  du^{i,N}_t &= dw^{i,N}_t + d(\grad q^{i,N}_t)\\
  &= \left[ -(u^N_t \cdot \grad) w^{i,N}_t + \nu \lap w^{i,N}_t - (\gradt u^N_t) w^{i,N}_t
  \right] \, dt - \sqrt{2\nu} \grad w^{i,N}_t dW^i_t\\ 
  &= \left[ -(u^N_t \cdot \grad) u^{i,N}_t + \nu \lap u^{i,N}_t - (\gradt u^N_t) u^{i,N}_t
  \right] \, dt - \sqrt{2\nu} \grad u^{i,N}_t dW^i_t +\\ 
  &\qquad +\left[ -(u^N_t \cdot \grad) (\grad q^{i,N}_t) + \nu \lap (\grad q^{i,N}_t)
    - (\gradt u^N_t) (\grad q^{i,N}_t) \right] \, dt -\\ 
  &\qquad - \sqrt{2\nu} \grad (\grad q^{i,N}_t) dW^i_t + d(\grad q^{i,N}_t). 
\end{align*}
If we define $P^{i,N}_t$ by
$$
P^{i,N}_t = \int_0^t \left[ (u^N_s \cdot \grad) q^{i,N}_s - \nu \lap
  q^{i,N}_s\right] \, ds + \int_0^t \sqrt{2\nu} \grad q^{i,N}_s \cdot dW^i_s +
q^{i,N}_t
$$
then we have
\begin{multline}\label{eDui}
  du^{i,N}_t + \left[ (u^N_t \cdot \grad) u^{i,N}_t - \nu \lap u^{i,N}_t + (\gradt u^N_t) u^{i,N}_t
  \right] \, dt \\+ d(\grad P^{i,N}_t) + \sqrt{2\nu} \grad u^{i,N}_t dW^i_t = 0. 
\end{multline}
Notice that $u^{i,N}$ is divergence free by definition, and thus  $\grad
u^{i,N} dW_i$ is also divergence free. Thus $d(\grad P^{i,N})$, the
only other term with possibly non-zero quadratic variation, must have a divergence free martingale part. Since the martingale part of $d(\grad P^{i,N})$ is also a gradient, it must be $0$. Thus
\begin{equation*}
  d(\grad P^{i,N}_t) = \grad p^{i,N}_t \, dt
\end{equation*}
for some function $p^{i,N}_t$, which proves \eqref{eSPDEuiN}.

The identity \eqref{eSPDEuN} now follows by summing \eqref{eSPDEuiN}
in $i$, dividing by $N$, and defining $p^N_t$ by
\begin{equation*}
  p^N_t =  \frac{1}{2} \grad \abs{u^N_t}^2 + \frac{1}{N} \sum_{i=1}^N p^{i,N}_t \qedhere
\end{equation*}
\end{proof}

\begin{proof}[Proof of Theorem \ref{tConvergence}]
  Let $u$ be a solution of the Navier-Stokes equations, with initial
  data $u_0$, and set $v^N = u^N - u$. Then $v^N$ satisfies the SPDE
  \begin{multline*}
    dv^N_t + (v^N_t \cdot \grad)u_t \, dt + (u_t \cdot \grad) v^N_t \, dt
    - \nu \lap v^N_t \, dt + \grad(p^N_t - p_t) \, dt +\\
    + \frac{\sqrt{2 \nu}}{N} \sum_{i=1}^N \grad u^{i,N}_t
    dW^i_t = 0,
  \end{multline*}
where $p$ is the pressure in the Navier-Stokes equations, and $p^N$ the pressure term in \eqref{eSPDEuN}.

Thus by It\^o's formula,
\begin{multline*}
  \frac{1}{2} d \lpnorm{v^N_t}{2}^2 + \ip{v^N_t}{(v^N_t \cdot \grad)
    u_t} \, dt + \nu \lpnorm{\grad v^N_t}{2}^2 \, dt +\\ 
  + \frac{2\nu}{N} \sum_{i=1}^N v^N_t \cdot (\grad u^{i,N}_t \, dW^{i,k}_t) =
  \frac{\nu}{N^2} \sum_{i=1}^N \lpnorm{\grad u^{i,N}_t}{2}^2 \, dt.
\end{multline*}
Here the notation $\ip{f}{g}$ denotes the $L^2$ innerproduct of $f$
and $g$. Taking expected values gives us
$$
\frac{1}{2} \del_t \E\lpnorm{v^N_t}{2}^2 \leq \frac{U}{L} \E
\lpnorm{v^N}{2}^2 + \frac{\nu U^2}{N L^2}
$$
and by Gronwall's lemma we have
$$
\E \lpnorm{v^N_t}{2}^2 \leq \frac{ 2 \nu U^2}{L^2 N} t e^{\frac{U t}{L}}
$$
concluding the proof.
\end{proof}

\section{Convergence as \texorpdfstring{$t \to \infty$}{t to infinity}}\label{sTtoInfty}
In this section, we fix $N$, and consider the behaviour of the system
\mceuler as $t \to \infty$. We show that the system \mceuler does not
dissipate all its energy as $t \to \infty$. Roughly speaking we show
$$
\limsup_{t \to \infty} \E \lpnorm{\grad u_t}{2}^2 \geq
O\left(\tfrac{1}{N}\right),
$$
with constants independent of viscosity.
This is in contrast to the true (unforced) Navier-Stokes equations,
which dissipate all of its energy as $t \to \infty$ (provided of
course the solutions are defined globally in time).

In general we are unable to compute exact asymptotic behaviour of the
system \mceuler as $t \to \infty$. But in the special case of shear
flows, we compute this exactly, and show that the system eventually
converges to a constant, retaining exactly $\frac{1}{N}$ times its
initial energy.

For the remainder of this section, we pick a fixed $N \in \N$ and for notational convenience we omit the superscript $N$. We begin by computing exactly the asymptotic behaviour of the system \mceuler in the special case of shear flows. 

\begin{proposition}\label{propLongTimeOmegaShear}
Suppose the initial data $u_0(x) = (\phi_0(x_2), 0)$ for some $\holderspace{1}{\alpha}$ periodic function $\phi$. If $u$ is the velocity field that solves the system \mceuler with initial data $u_0$, then
\begin{equation}\label{eqnLongTimeOmegaShear}
\lim_{t \to \infty} \E\,\omega_t(x)^2 = \frac{1}{N} \lpnorm{\omega_0}{2}^2
\end{equation}
where $\omega = \curl u$ is the vorticity.
\end{proposition}
\begin{proof}
Let $X^i$, $Y^i$ be the flows in the system \mceuler, and as before define $u^i$ to be the $i^\text{th}$ summand in \eqref{euN}, and $\omega^i = \omega_0 \circ Y^i$.

First note that the SPDE's for $u^i$ and $u$ (equations \eqref{eSPDEuiN} and \eqref{eSPDEuN}) are all translation invariant. Thus since the initial data is independent of $x_1$, the same must be true for all time. Since $u^i$, $u$ are divergence free, the second coordinate must necessarily be $0$, and the form of the initial data is preserved. Namely,
$$
u_t(x) = (\phi_t(x_2), 0) \quad\text{and}\quad u^i_t(x) = (\phi^i_t(x_2), 0) 
$$
for some $C([0,\infty), \holderspace{1}{\alpha})$ periodic functions $\phi^i, \phi$.

Now the definition of $X^i$ shows that
$$
X^i_t(y) =
\begin{pmatrix}
  y_1 + \lambda^{i,1}_t(y_2)\\
  y_2 + \sqrt{2\nu} W^{i,2}_t
\end{pmatrix}
$$
and hence
$$
Y^i_t(x) =
\begin{pmatrix}
  x_1 + \mu^{i,1}_t(x_2)\\
  x_2 - \sqrt{2\nu} W^{i,2}_t
\end{pmatrix}.
$$
Recall $\lambda^i$, $\mu^i$ are as in \eqref{eLambdai}, \eqref{eMui}, and here the notation $\lambda^{i,1}$ to denotes the first coordinate of $\lambda^i$. This immediately shows
\begin{align*}
\omega^i_t(x) &= \omega_0 \circ Y^i_t (x)\\
    &= -\del_2 \phi_0( x_2 - \sqrt{2\nu} W^{i,2}_t )
\end{align*}
where $W^{i,2}$ again denotes the second coordinate of the Brownian motion $W^i$.

Now using standard mixing properties of Brownian motion~\cite{kunita}*{Section 1.3}, (or explicitly computing in this case) we know that for every $x \in [0,L]^2$
\begin{equation}\label{eqnEOmegaSq}
  \lim_{t\to\infty} \E\, \omega^i_t(x)^2 = \lim_{t\to\infty} \E\left[-\del_2 \phi_0( x_2 - \sqrt{2\nu} W^{i,2}_t )\right]^2 = \lpnorm{\omega_0}{2}^2
\end{equation}
and
\begin{align}
\nonumber
\lim_{t\to\infty} \E\, \omega^i_t(x) \omega^j_t(x) &= \lim_{t\to\infty} \E \left[\del_2 \phi_0( x_2 - \sqrt{2\nu} W^{i,2}_t ) \del_2 \phi_0( x_2 - \sqrt{2\nu} W^{j,2}_t )\right]\\
\nonumber
    &= \lim_{t\to\infty} \left[\E \del_2 \phi_0( x_2 - \sqrt{2\nu} W^{i,2}_t )\right] \, \left[\E \del_2 \phi_0( x_2 - \sqrt{2\nu} W^{j,2}_t )\right]\\
\nonumber
    &= \left(\frac{1}{L^2} \int_{[0,L]^2} \del_2 \phi_0\right)^2\\
\label{eqnEOmegaiOmegaj}  &= 0
\end{align}
when $i \neq j$.

Now by two dimensional Cauchy formula \eqref{eCauchy2D}
$$
\omega_t = \frac{1}{N} \sum_{i=1}^N \omega_0 \circ Y^i_t.
$$
(since in our case, $\omega^1_0 = \cdots = \omega^N_0 = \omega_0$). Thus
\begin{equation}
  \E\,\omega_t(x)^2 = \frac{1}{N^2} \sum_{i=1}^N \E\,\omega_0 \circ
  Y^i_t (x)^2 + \frac{2}{N^2} \sum_{i=2}^N \sum_{j = 1}^{i-1} \E\,
  \omega_0 \circ Y^i_t(x) \omega_0 \circ Y^j_t(x) 
\end{equation}
and using \eqref{eqnEOmegaSq} and \eqref{eqnEOmegaiOmegaj} the proof is complete.
\end{proof}

We remark that all we need for \eqref{eqnLongTimeOmegaShear} to hold is the identities \eqref{eqnEOmegaSq} and \eqref{eqnEOmegaiOmegaj}. Equality \eqref{eqnEOmegaSq} is guaranteed provided reasonable ergodic properties of the flow $X^i_t$ are known. Equality \eqref{eqnEOmegaiOmegaj} is guaranteed provided the flows $X^i_t$ and $X^j_t$ eventually decorrelate.

While we are unable to guarantee these properties for a more general class of flows, we conclude this section by proving a weaker version of \eqref{eqnLongTimeOmegaShear} for two dimensional flows with general initial data.

\begin{theorem}\label{tLimSupGradULowerBound}
  Let $X^i, u$ be a solution to the system \mceuler with (spatial)
  mean zero initial data $u_0 \in \holderspace{1}{\alpha}$ and
  periodic boundary conditions. Suppose further $u \in C( [0,\infty),
  \holderspace{1}{\alpha} )$. Then
\begin{equation}\label{eLimSupGradULowerBound}
\limsup_{t \to \infty} \E \lpnorm{\grad u_t}{2}^2 \geq \frac{1}{N L^2}
\lpnorm{u_0}{2}^2 
\end{equation}
\end{theorem}

Note that the assumption $u \in C( [0,\infty),
\holderspace{1}{\alpha})$ is satisfied in the two dimensional
situation with $\holderspace{1}{\alpha}$ initial data (Theorem
\ref{t2DGlobalExistence}). The proof we provide below will also work
in the three dimensional situation, as long as global existence and
well-posedness of \mceuler is known.

As is standard with the Navier-Stokes equations, the condition that
$u_0$ is (spatially) mean zero is not a restriction. By changing
coordinates to a frame moving with the mean of the initial velocity,
we can arrange that the initial data (in the new frame) has spatial
mean $0$.

Finally we remark that the lower bound in inequality \eqref{eLimSupGradULowerBound} is sharp, since in the special case of shear flows we have the \textit{equality} \eqref{eqnLongTimeOmegaShear}. However we are unable to obtain a bound on $\liminf \E \lpnorm{\grad u_t}{2}^2$.

\begin{proof}[Proof of Theorem \ref{tLimSupGradULowerBound}]
  As before let $u^i_t = \lhp\left[ (\gradt Y^i_t) u_0 \circ Y^i_t
  \right]$. Using Lemma \ref{lSPDEui} and It\^o's formula we have
\begin{multline}\label{eDuiSquared}
  \frac{1}{2} d \abs{u^i}^2 + u^i \cdot \left[ (u \cdot \grad) u^i -
    \nu \lap u^i + (\gradt u) u^i + \grad p^i \right] \, dt +\\ 
  + \sqrt{2 \nu} u^i \cdot (\grad u^i\, dW^i) = \nu \abs{\grad u^i}^2
  \, dt
\end{multline}
Note that
$$
\int u^i \cdot ((\gradt u) u^i) = \int ((\grad u) u^i) \cdot u^i =
\int u^i \cdot (u \cdot \grad) u^i = 0 
$$
Thus integrating \eqref{eDuiSquared} in space, and using the fact that
$u^i$ is divergence free gives
$$
d \lpnorm{u^i}{2}^2 = 0
$$
and hence $\lpnorm{u^i_t}{2} = \lpnorm{u^i_0}{2}$ almost
surely.\smallskip

Now suppose that for some $\epsilon > 0$, there exists $t_0$ such that for all $t > t_0$
$$
\E \lpnorm{ \grad u_t}{2}^2 \leq \frac{1}{N L^2} \lpnorm{u_0}{2}^2 - \epsilon.
$$
Using It\^o's formula and \eqref{eSPDEuN} gives
\begin{multline*}
  \frac{1}{2} d \abs{u}^2 + u \cdot \left[ (u \cdot \grad) u - \nu
    \lap u + \grad p \right] \, dt 
  + \frac{\sqrt{2 \nu}}{N} \sum_{i=1}^N u \cdot (\grad u^i\, dW^i) \\=
  \frac{\nu}{N^2} \sum_{i=1}^N \abs{\grad u^i}^2 \, dt
\end{multline*}
Integrating in space, and taking expected values gives
\begin{align*}
  \frac{1}{2\nu} \del_t \E \lpnorm{u_t}{2}^2 &= \E \left[
    \frac{1}{N^2} \sum_{i=1}^N \lpnorm{\grad u^i_t}{2}^2 -
    \lpnorm{\grad u_t}{2}^2 \right]\\ 
  &\geq \E \left[ \frac{1}{N^2 L^2} \sum_{i=1}^N \lpnorm{u^i_t}{2}^2 - \lpnorm{\grad u_t}{2}^2 \right]\\
  &= \frac{1}{N^2 L^2} \sum_{i=1}^N \lpnorm{u_0}{2}^2 - \E \lpnorm{\grad u_t}{2}^2\\
  &= \frac{1}{N L^2} \lpnorm{u_0}{2}^2 - \E \lpnorm{\grad u_t}{2}^2\\
  &\geq \epsilon
\end{align*}
for $t \geq t_0$. Here we used the Poincar\'e inequality to obtain the
second inequality above. Note that we have assumed that the initial
data has (spatial) mean $0$. Since the (spatial) mean is conserved by
the system \mceuler, $u_t$ also has (spatial) mean zero, and our
application of the Poincar\'e inequality is valid.

Now, the above inequality immediately implies $\E \lpnorm{ u_t }{2}^2$
becomes arbitrarily large as $t \to \infty$. This is a contradiction
because
$$
\lpnorm{u_t}{2} = \lpnorm{\frac{1}{N} \sum_{i=1}^N u^i_t }{2} \leq \frac{1}{N}
\sum_{i=1}^N \lpnorm{u^i_t}{2} = \lpnorm{u_0}{2}
$$
holds almost surely.
\end{proof}
\appendix
\section{Local existence.}\label{aLocalExistence}
In this appendix we provide the proof of Proposition
\ref{pLocalExistence}. A similar proof appeared in \cite{sperturb}
(see also \cite{ele}), and the proof provided here is based on similar
ideas. We present the proof here because we require local existence
for $\holderspace{1}{\alpha}$ initial data (the proof in
\cite{sperturb} used $\holderspace{2}{\alpha}$), and to ensure that
bounds and existence time therein are independent of $N$.

Without loss, we assume $t_0 = 0$, and $u_0^1$, \dots, $u_0^N$ to be
$\F_0$ measurable. Let $U$ be a large constant and $T$ a small time,
both of which will be specified later.

Define $\U = \U(T,U)$ be the set of all time continuous $\F_t$-adapted
$\holderspace{1}{\alpha}$ valued divergence free and spatially
periodic processes $u$ such that
$$
u_0 = \frac{1}{N} \sum_{i=1}^N u^i_0 \quad\text{and}\quad \sup_{t \in
  [0,T]} \holdernorm{u_t}{1}{\alpha} \leq U
$$
hold almost surely. Also, we define $\M = \M(T)$ to be the set of all
time continuous $\F_t$-adapted $\holderspace{1}{\alpha}$ valued
spatially periodic processes $\mu$ such that
$$
\mu_0 = 0 \quad\text{and}\quad \sup_{t \in [0,T]} \holdernorm{\grad
  \mu_t}{0}{\alpha} \leq \frac{1}{2}
$$
hold almost surely.

Now given $u \in \U$ let $X^{i,u}$ be the flow solving the SDE
$$
dX^{i,u}_t = u_t(X^{i,u}_t) \, dt + \sqrt{2\nu} \, dW^i_t
$$
with initial data $X^{i,u}_0(y) = y$. As before, define $Y^{i,u}_t =
(X^{i,u}_t)\inv$, and define
\begin{gather*}
\lambda^{i,u}_t(y) = X^{i,u}_t(y) - y\\
\mu^{i,u}_t(x) = Y^{i,u}_t(x) - x
\end{gather*}
to be the Eulerian and Lagrangian displacements respectively.

Finally define the (non-linear) operator $\W$ by
$$
\W(u)_t = \frac{1}{N} \sum_{i=1}^N \lhp\left[ \left( \gradt Y^{i,u}_t
  \right)\, \left( u^i_0 \circ Y^{i,u}_t \right) \right] 
$$
Clearly a fixed point of $\W$ will produce a solution to the system
\mceulert. Thus the proof will be complete if we show that for an
appropriate choice of $T$ and $U$, $\W$ maps $\U$ into itself, and is a
contraction with respect to the weaker norm
$$
\norm{u}_\U = \sup_{t \in [0,T]} \holdernorm{u_t}{0}{\alpha}
$$

We first show $\W$ maps $\U$ into itself, using the two lemmas below.
\begin{lemma}\label{lWbound}
There exists $c = c(\alpha)$ such that
$$
\holdernorm{\W(u)}{1}{\alpha} \leq c \left[ \max_{1 \leq i \leq N}
  \left(1 + \holdernorm{\grad \mu^{i,u}}{0}{\alpha}\right)^{2+\alpha}
\right] \frac{1}{N} \sum_{i=1}^N \holdernorm{u^i_0}{1}{\alpha}
\quad\as
$$
\end{lemma}
\begin{proof}
First recall $\lhp$ vanishes on gradients. Thus
\begin{equation}\label{eIntByParts}
\lhp\left[ (\gradt Y) \, v \right] = - \lhp\left[ (\gradt v) \, Y \right].
\end{equation}
Now
\begin{align*}
  \del_i \lhp\left[ (\gradt Y)v\right] &= \lhp\left[ (\gradt Y) \del_i v + (\gradt \del_i Y) v \right]\\
  &= \lhp\left[ (\gradt Y) \del_i v - (\gradt v) \del_i Y \right]
\end{align*}
where we used \eqref{eIntByParts} for the second term. Note that the
right hand side involves only first order derivatives. Since $\lhp$ is
a standard Calder\'on-Zygmund singular integral operator, which is
bounded on H\"older spaces, we obtain the estimate
$$
\holdernorm{ \lhp\left[ (\gradt Y) v \right] }{1}{\alpha} \leq c
\holdernorm{\gradt Y}{0}{\alpha} \holdernorm{v}{1}{\alpha}
$$
for some constant $c = c(\alpha)$.

Applying this estimate to $\W$, we have
\begin{equation}\label{eWu1}
  \holdernorm{\W(u)_t}{1}{\alpha} \leq \frac{c}{N} \sum_{i=1}^N
  \holdernorm{\gradt Y^{i,u}_t}{0}{\alpha} \holdernorm{ u^i_0 \circ
    Y^{i,u}_t }{1}{\alpha}\quad\as 
\end{equation}
from which the Lemma follows.
\end{proof}
\begin{lemma}\label{lLambdaMuBound}
There exists $T = T(U, \alpha)$ such that $\lambda^{i,u}, \mu^{i,u} \in \M(T)$.
\end{lemma}

We note that the diffusion coefficient is spatially constant, and thus
we get the desired (almost sure) control on $\grad \lambda$. Since
$\divergence u = 0$, $\det(\grad X^{i,u}) = 1$, giving the desired
control on $\grad \mu$. The details are standard, and we do not
provide the them here (see for instance
\cites{sperturb,kunita}).\medskip

Now choosing $U = c(\frac{3}{2})^{2+\alpha} \frac{1}{N} \sum
\holdernorm{u^i_0}{1}{\alpha}$, and then choosing $T$ as in Lemma
\ref{lLambdaMuBound}, Lemma \ref{lWbound} shows that $\W$ maps
$\U(U,T)$ into itself. Note that each summand on the right of
\eqref{eWu1} is bounded by $U$, which will prove the bound
\eqref{eUiBound}. Further, given a uniform (in $i$) bound on
$\holdernorm{u^i_0}{1}{\alpha}$, our choice of $U$ can be made
independent of $N$.

It remains to show that $\W$ is a contraction. By definition of $\W$ we have
\begin{align*}
  \W(u)_t - \W(v)_t &= \frac{1}{N} \sum_{i=1}^N \lhp\left[ (\gradt Y^{i,u}_t)
    u^i_0 \circ Y^{i,u}_t - (\gradt Y^{i,v}_t) u^i_0 \circ Y^{i,v}_t
  \right]\\ 
  &= \frac{1}{N} \sum_{i=1}^N \lhp\left[ (\gradt Y^{i,u}_t) \left( u^i_0 \circ
      Y^{i,u}_t  - u^i_0 \circ Y^{i,v}_t \right)\right] +\\ 
  &\qquad + \frac{1}{N}\sum_{i=1}^N \lhp\left[ \left( \gradt Y^{i,u}_t -
      \gradt Y^{i,v}_t \right) u^i_0 \circ Y^{i,v}_t \right]\\ 
  &= \frac{1}{N} \sum_{i=1}^N \lhp\left[ (\gradt Y^{i,u}_t) \left( u^i_0 \circ
      Y^{i,u}_t  - u^i_0 \circ Y^{i,v}_t \right)\right] -\\ 
  &\qquad- \frac{1}{N}\sum_{i=1}^N \lhp\left[ \gradt (u^i_0 \circ Y^{i,v}_t)
    \left( Y^{i,u}_t - Y^{i,v}_t \right) \right]
\end{align*}
where we used the identity \eqref{eIntByParts} to obtain the last
equality. Now we recall that $\mu^{i,u}, \mu^{i,v} \in \M$, and take
$\holderspace{0}{\alpha}$ norms. This gives
\begin{equation}\label{eWuMinusWv1}
  \holdernorm{\W(u)_t - \W(v)_t}{0}{\alpha} \leq \frac{c}{L N} \sum_{i=1}^N
  \holdernorm{u^i_0}{1}{\alpha} \holdernorm{Y^{i,u}_t -
    Y^{i,v}_t}{0}{\alpha} \quad\as 
\end{equation}

Now from the definition of $Y^{i,u}$ and $Y^{i,v}$ we have
$$
Y^{i,u}_t - Y^{i,v}_t = \int_0^t \left[ u_s(Y^{i,u}_s) -
  v_s(Y^{i,v}_s) \right] \, ds. 
$$
Taking $\holderspace{0}{\alpha}$ norms, and applying Gronwall's
inequality, and absorbing the exponential in time factor into the
constant $c$ gives
$$
\holdernorm{Y^{i,u}_t - Y^{i,v}_t}{0}{\alpha} \leq c \int_0^t
\holdernorm{u_s - v_s}{0}{\alpha} \, ds \quad\as
$$

Returning to \eqref{eWuMinusWv1} we have
$$
\holdernorm{ \W(u)_t - \W(v)_t }{0}{\alpha} \leq \frac{c t}{L}
\sup_{s\leq t} \holdernorm{u_s - v_s}{0}{\alpha} \frac{1}{N} \sum_{i=1}^N
\holdernorm{ u^i_0 }{1}{\alpha} \quad\as
$$
Choosing $t$ small enough one can ensure $\W$ is a contraction
mapping. A standard iteration now shows the existence of a fixed point
of $\W$, concluding the proof of Proposition \ref{pLocalExistence}.

\section{Logarithmic \texorpdfstring{$L^\infty$}{L-infinity} bound on
  singular integral operators.}\label{aSIOlogBound}
In this appendix we provide a proof of Lemma \ref{lSIOlogBound}. We
restate it here for the readers convenience.
\begin{lemma*}[\ref{lSIOlogBound}]
  If $u$ is divergence free and periodic in $\R^d$, then for any
  $\alpha \in (0,1)$, there exists a constant $c = c(\alpha, d)$ such
  that
$$
\linfnorm{\grad u} \leq c \linfnorm{\omega} \left( 1 + \ln^+ \left(
    \frac{ \holdersnorm{\omega}{\alpha} }{ \linfnorm{\omega} } \right)
\right) 
$$
where $\omega = \curl u$.
\end{lemma*}
\begin{proof}
  Let $K$ be a standard Calder\'on-Zygmund periodic kernel, and we
  define the operator $T$ by
$$ Tf = K * f $$
then we will prove that
\begin{equation}\label{eLinfBoundf}
\linfnorm{T f} \leq c \linfnorm{f} \left( 1 + \ln^+ \left( \frac{ \holdersnorm{f}{\alpha} }{ \linfnorm{f} } \right) \right)
\end{equation}
This immediately implies Lemma \ref{lSIOlogBound} because we know
$\grad u = -\grad (\lap\inv) \curl \omega$, and $-\grad (\lap\inv)
\curl$ is a Calder\'on-Zygmund type singular integral operator.

Now we prove \eqref{eLinfBoundf}. We assume for convenience that all
functions are periodic on the cube $[0,1]^d$. We also recall that the
kernel $K$ satisfies the properties
\begin{enumerate}
\item\label{aKdecay} $K(y) \leq c \abs{y}^{-d}$ when $\abs{y} \leq \frac{1}{2}$.
\item\label{aKcancellation} $\int_{\abs y = r} K(y) \, d\sigma(y) = 0$
  for any $r \in (0,\frac{1}{2})$.
\end{enumerate}

Pick any $\epsilon \in (0,\frac{1}{2})$. Then
\begin{align}
\nonumber Tf(x) &= \int_{[0,1]^d} K(y) f(x-y) \, dy\\
\label{eTfsplit}   &\leq \int_{\abs y < \epsilon} K(y) f(x-y) \, dy + \int_{\abs y \geq \epsilon} K(y) f(x-y) \, dy
\end{align}
Using property \ref{aKdecay} about $K$, we bound the second integral
by
\begin{align*}
  \abs{\int_{\abs y \geq \epsilon} K(y) f(x-y) \, dy} &\leq c \linfnorm{f} \int_{r=\epsilon}^1 \frac{1}{r^d} r^{d-1} \, dr\\
  &\leq c \linfnorm{f} \ln \left( \tfrac{1}{\epsilon} \right)
\end{align*}
Using property \ref{aKcancellation} about $K$, we bound the first integral in \eqref{eTfsplit} by
\begin{align*}
  \abs{\int_{\abs y < \epsilon} K(y) f(x-y) \, dy} & = \abs{\int_{\abs y < \epsilon} K(y) (f(x-y) - f(x)) \, dy}\\
  &\leq c \holdersnorm{f}{\alpha} \int_{\abs y < \epsilon} \abs{y}^{\alpha - d} \, dy\\
  &= c \holdersnorm{f}{\alpha} \epsilon^\alpha
\end{align*}
Combining estimates we have
$$
\linfnorm{Tf} \leq c\left[ \epsilon^\alpha \holdersnorm{f}{\alpha} +
  \linfnorm{f} \ln \left( \tfrac{1}{\epsilon} \right) \right]
$$
Choosing
$$
\epsilon = \min\left\{ \frac{1}{2},
  \left(\frac{\linfnorm{f}}{\holdersnorm{f}{\alpha}}\right)^{1/\alpha}
\right\}
$$
finishes the proof.
\end{proof}

\end{document}